# A Quick Introduction To Dwork's Conjecture

Daqing Wan

ABSTRACT. This paper is an expanded version of my lecture delivered at the 97 Seattle summer research conference on finite fields. It gives a quick exposition of Dwork's conjecture about $p$-adic meromorphic continuation of his unit root L-function arising from a family of algebraic varieties defined over a finite field of characteristic $p$. As a simple illustration, we discuss the classical example of the universal family of elliptic curves where the conjecture is already known to be true and where the conjecture is closely related to arithmetic of modular forms such as the Gouvêa-Mazur conjecture. Special attention is given to questions related to the $p$-adic absolute values of the unit root L-function. In particular, it is observed that an average version of a suitable $p$-adic Riemann hypothesis is true for the elliptic family. Following a suggestion of Mike Fried, I also include a section describing some of my personal interactions with Dwork. This extra section serves as a dedication to the memory of Dwork who actively attended the conference and died nine months later.

## 1. One version of Dwork's conjecture

This section gives a quick reformulation of Dwork's conjecture in the general case. There are several different but essentially equivalent languages, such as $p$-adic Galois representations, $p$-adic ètale sheaves and unit root F-crystals, that could be used to describe Dwork's conjecture. To be compatible with the general theme of the conference, I will use the language of $p$-adic representations and Galois groups. This provides a short although not the simplest reformulation of Dwork's conjecture. In one sentence, the conjecture simply says that if $\rho$ is a continuous $p$-adic Galois representation coming from algebraic geometry over a finite field of characteristic $p$, then the L-function $L(\rho, T)$ is $p$-adic meromorphic. We now make this a little more precise.

Let $q$ be a power of a fixed prime number $p$ and let $\mathbb{F}_q$ be the finite field of $q$ elements. For a geometrically connected algebraic variety $X$ defined over $\mathbb{F}_q$, let $\pi_1^{\text{arith}}(X)$ denote the arithmetic fundamental group of $X$. This is a profinite group. More precisely, $\pi_1^{\text{arith}}(X)$ is the profinite completion of the finite Galois groups of pointed finite unramified Galois coverings of $X$. If $X$ is integral and normal with function field $\mathbb{F}_q(X)$, then $\pi_1^{\text{arith}}(X)$ is simply the profinite Galois group $\text{Gal}(\mathbb{F}_q(X)^{\text{sep}}/\mathbb{F}_q(X))$ modulo the normal subgroup generated by the inertia subgroups $I_x$ at the closed points $x$ of $X$, where $\mathbb{F}_q(X)^{\text{sep}}$ denotes a fixed separable closure of $\mathbb{F}_q(X)$. We shall be interested in arithmetic and analytic properties of the L-function attached to a continuous $p$-adic representation of $\pi_1^{\text{arith}}(X)$.

Let $R$ be the ring of integers in a finite extension of the $p$-adic rational numbers $\mathbb{Q}_p$. For example, one could take $R = \mathbb{Z}_p$, the $p$-adic rational integers. Let $\rho$ be a continuous $p$-adic representation

$$\rho: \pi_1^{\text{arith}}(X) \longrightarrow \text{GL}_n(R).$$

The L-function of $\rho$ is defined as usual by the following infinite Euler product

$$L(\rho, T) = \prod_{x \in X_0} \frac{1}{\det(I - T^{\deg(x)}\rho(\text{Frob}_x))}, \tag{1.1}$$

where $X_0$ denotes the set of closed points of $X$ over $\mathbb{F}_q$ and $\text{Frob}_x$ denotes the Frobenius conjugacy class at $x$ of $\pi_1^{\text{arith}}(X)$. It is clear that the L-function $L(\rho, T)$ is a formal power series with coefficients in $R$. Thus, it is trivially $p$-adic analytic in the open unit disk $|T|_p < 1$.

For arithmetic applications, it is important to understand the $p$-adic meromorphic continuation of $L(\rho, T)$ and the nature of its zeroes. For instance, if $\rho$ is of finite order, the general theorem of Dwork-Grothendieck shows that $L(\rho, T)$ is a rational function. In this case, the zeroes and poles of $L(\rho, T)$ are integral powers of $\sqrt{q}$ by Deligne's theorem [De2] on Riemann hypothesis over finite fields. The total number of zeroes and poles of $L(\rho, T)$ can be bounded explicitly by Bombieri-Sperber [BS]. The full potential of the available theories has, however, not been fully exploited. In particular, there is still a great deal of work to be done to understand the $p$-adic absolute values of the zeroes and poles of $L(\rho, T)$, see Mazur [Ma], Oort [Oo], Illusie [Il] and [W1] for some positive results in this direction.

For infinite order $p$-adic representation $\rho$, the situation is naturally more complicated. First, the L-function $L(\rho, T)$ will not be rational in general unless $\rho$ is of very special type. The Dwork-Monsky [Mo] trace formula implies that $L(\rho, T)$ is $p$-adic meromorphic if $\rho$ is overconvergent in some sense. In terms of F-crystals, our overconvergent condition on $\rho$ simply means that the Frobenius map of the F-crystal is overconvergent with respect to some lifting. This condition is much weaker than Berthelot's overconvergent F-crystal [Be] which assumes that both the Frobenius map and the horizontal connection are overconvergent. All finite order representations are overconvergent in Berthelot's sense. It would be interesting to give a representation theoretic and/or group theoretic characterization of our overconvergent condition purely in terms of the representation $\rho$ itself without using its F-crystal counterpart.

Using the Monsky trace formula and a simple limiting argument, it is not hard to show that the L-function $L(\rho, T)$ has a $p$-adic meromorphic continuation to the closed unit disk $|T|_p \leq 1$. One plausible conjecture of Katz [K1] says that the unit (slope zero) part of $L(\rho, T)$ is given by the characteristic "polynomial" of the geometric Frobenius acting on the compact $p$-adic ètale cohomology of the $p$-adic ètale sheaf $\rho$ on $X$. This is known to be true if $\rho$ is the trivial character, see [ES]. Thus, at least conjecturally, the slope zero part of $L(\rho, T)$ is well understood. The higher slope portion of $L(\rho, T)$ is more difficult. In this direction, a general conjecture of Katz [K1] says that $L(\rho, T)$ is $p$-adic meromorphic everywhere. This conjecture is disproved in [W2]. In studying $p$-adic analytic variation [D5-6] of a family of varieties over $\mathbb{F}_q$, one naturally leads to $p$-adic representation $\rho$ which is not overconvergent in general and thus goes beyond the classical overconvergent theory. Dwork's conjecture says that the L-function of such a $p$-adic representation coming from algebraic geometry over a finite field of characteristic $p$ is still $p$-adic meromorphic.

For simplicity of description, we shall restrict our attention to an essential case. A continuous $p$-adic representation $\rho$ of $\pi_1^{\mathrm{arith}}(X/\mathbb{F}_q)$ is called **geometric** if $\rho$ arises as a relative $p$-adic ètale cohomology $R^i f_* \mathbb{Z}_p$ of a family of smooth proper varieties $f : Y \to X$ defined over $\mathbb{F}_q$. One could extend the notion of geometric representations in various ways. For instance, one could allow singular open families. One could also replace the constant $p$-adic sheaf $\mathbb{Z}_p$ on $Y$ by a more general finite order $p$-adic sheaf on $Y$. One could consider the tensor category generated by geometric $p$-adic representations. One version of Dwork's conjecture can be reformulated as follows.

CONJECTURE 1.1 (DWORK [D6]). *Let $\rho$ be a geometric $p$-adic representation of $\pi_1^{arith}(X)$. Then, the L-function $L(\rho, T)$ is $p$-adic meromorphic.*

Note that the $\ell$-adic analogue of Conjecture 1.1 is already known to be true if $\ell$ is a prime number different from the characteristic $p$. In fact, Grothendieck's rationality theorem [Gr] says that the L-function of any continuous $\ell$-adic representation of $\pi_1^{\mathrm{arith}}(X/\mathbb{F}_q)$ is always a rational function, whether the $\ell$-adic representation is geometric or not. The situation in the $p$-adic case is much more complicated. The intuition is that there are too many $p$-adic representations (generally quite wildly ramified at infinities) while there are relatively few $\ell$-adic representations (almost tame). But wildly ramified $p$-adic representations do arise from algebraic geometry and are important for arithmetic applications as we shall see in next section. One arithmetic implication of Conjecture 1.1 is the existence of a general $p$-adic equi-distribution theorem, which means that the zeroes with a given slope of the zeta function of a variety over a finite field are equi-distributed in a suitable $p$-adic sense when the variety moves through an algebraic family.

As mentioned above, Conjecture 1.1 is known to be true when $\rho$ is overconvergent. This is the only case for which Conjecture 1.1 has been proved. Unfortunately, such overconvergent geometric representations are at present quite rare. The known cases include the universal family of ordinary plane curves of genus $g$ with $g \leq 3$ and a certain family of ordinary K-3 surfaces [D6]. Under the much stronger condition that $\rho$ is a unit root overconvergent F-crystal in Berthelot's sense, the L-function $L(\rho, T)$ is known to be a rational function by recent work of de Jong-Berthelot-Tsuzuki [Dj] on finite dimensionality of rigid cohomology. In this case, the unit root overconvergent F-crystal $\rho$ extends, after a finite covering $X' \to X$, to a projective completion of $X'$ and thus $\rho$ is almost tame as in $\ell$-adic case. In a series of future articles, we will prove Conjecture 1.1 in full generality, see [W4-5] for a complete proof in an essential non-overconvergent setup.

Once we know that the geometric L-function $L(\rho, T)$ is $p$-adic meromorphic. Many fundamental questions arise, such as order of poles, special values, and most importantly the $p$-adic absolute values of its zeroes. In next section, we treat the simplest elliptic family case of Conjecture 1.1 and discuss its relation to the Gouvêa-Mazur conjecture about dimension variation of $p$-adic modular forms. Our emphasis will be on various attempts to understand the $p$-adic absolute values of the zeroes and poles, namely, some sort of $p$-adic Riemann hypothesis or $p$-adic Ramanujan-Peterson conjecture. These questions are not well understood at present. The depth of the elliptic example should give an indication about the potential significance of Dwork's conjecture in the general case. It also suggests a possible general connection with arithmetic of automorphic forms. In particular, the Kloosterman family case of Conjecture 1.1 should be related to arithmetic of Maass forms.

## 2. The elliptic family and modular forms

In this section, we let $\mathbb{F}_q$ be the prime field $\mathbb{F}_p$ and assume that $p > 2$. Consider the Legendre family $E_x$ of elliptic curves whose affine equation is given by

$$E_x : y_1^2 = y_2(y_2 - 1)(y_2 - x),$$

where $x \in \mathbb{A}^1 - \{0, 1\}$. This is the universal elliptic curve of level 2 parametrized by $x \in \mathbb{A}^1 - \{0, 1\}$. We explain explicitly Dwork's conjecture for this family and its relation to $p$-adic modular forms. In terms of exposition style, we use the simplest formula approach to be as self-contained as possible.

For each geometric point $x \in \mathbb{F}_{p^{\deg(x)}} - \{0, 1\}$, the fibre $E_x$ is an elliptic curve defined over the finite field $\mathbb{F}_{p^{\deg(x)}}$. The zeta function of $E_x$ over the finite field $\mathbb{F}_{p^{\deg(x)}}$ is defined by the infinite product

$$Z(E_x, T) = \prod_{y \in (E_x)_0} \frac{1}{1 - T^{\deg(y)}},$$

where $(E_x)_0$ denotes the set of closed points of $E_x/\mathbb{F}_{p^{\deg(x)}}$. Alternatively, the zeta function can be defined as a generating function for the number of rational points over various extension fields of $\mathbb{F}_{p^{\deg(x)}}$:

$$Z(E_x, T) = \exp\left(\sum_{k=1}^{\infty} \frac{\#E_x(\mathbb{F}_{p^{k\deg(x)}})}{k} T^k\right).$$

It is well known that $Z(E_x, T)$ is a rational function of the following form

$$Z(E_x, T) = \frac{P_x(T)}{(1-T)(1-p^{\deg(x)}T)},$$

where $P_x(T)$ is a quadratic polynomial with coefficients in $\mathbb{Z}$. If we factor $P_x(T)$ over the algebraic closure $\bar{\mathbb{Q}}$ of $\mathbb{Q}$, we can write

$$P_x(T) = (1 - \alpha(x)T)(1 - \beta(x)T),$$

where $\alpha(x)$ and $\beta(x)$ are algebraic integers. The functional equation shows that

$$\alpha(x)\beta(x) = p^{\deg(x)}.$$

Thus, in order to have a complete understanding of the zeta function $Z(E_x, T)$, it suffices to understand one of the reciprocal roots of $P_x(T)$, say, $\alpha(x)$. The functional equation also implies that for each prime $\ell \neq p$, both $\alpha(x)$ and $\beta(x)$ are $\ell$-adic unit. The complex absolute value is given by Hasse's theorem on Riemann hypothesis:

$$|\alpha(x)| = |\beta(x)| = \sqrt{p}^{\deg(x)}.$$

The $p$-adic absolute values of the roots are a little more complicated to describe.

Define the Hasse polynomial $H(x)$ by

$$H(x) = \sum_{i=0}^{(p-1)/2} \binom{(p-1)/2}{i}^2 x^i.$$

This is a polynomial of degree $(p-1)/2$ with distinct roots over the algebraic closure of $\mathbb{F}_p$. We fix an embedding of $\bar{\mathbb{Q}}$ into the completion $\Omega_p$ of an algebraic closure of $\mathbb{Q}_p$. Arrange the reciprocal roots $\alpha(x)$ and $\beta(x)$ such that

$$0 \leq \mathrm{ord}_p \alpha(x) \leq \mathrm{ord}_p \beta(x).$$

In the supersingular case that $H(x)$ vanishes at $x$ (there are only $(p-1)/2$ such $x$), we have

$$\mathrm{ord}_p \alpha(x) = \mathrm{ord}_p \beta(x) = \frac{\deg(x)}{2}.$$

In the ordinary case that $H(x)$ does not vanish at $x$, then we have

$$\mathrm{ord}_p \alpha(x) = 0, \ \mathrm{ord}_p \beta(x) = \deg(x).$$

For ordinary $E_x$, the first $p$-adic ètale cohomology group $H^1_{\mathrm{et}}(E_x, \mathbb{Z}_p)$ is isomorphic to $\mathbb{Z}_p$ (the dual of the $p$-adic Tate module). The eigenvalue of the geometric Frobenius $\mathrm{Frob}_x$ acting on $H^1_{\mathrm{et}}(E_x, \mathbb{Z}_p)$ is precisely given by the $p$-adic unit $\alpha(x)$.

Let $X$ be the punctured affine line $\mathbb{A}^1 - \{0, 1, H(x) = 0\}$ defined over $\mathbb{F}_p$. The family

$$f : E_x \longrightarrow x \in X$$

of elliptic curves $E_x$ parametrized by $x \in X$ is the universal family of ordinary elliptic curves of level 2 over $\mathbb{F}_p$. The first relative $p$-adic ètale cohomology $R^1 f_* \mathbb{Z}_p$ of the family $f$ gives a continuous rank one $p$-adic representation of $\pi_1^{\mathrm{arith}}(X)$:

$$\rho_E : \ \pi_1^{\mathrm{arith}}(X) \longrightarrow \mathrm{GL}_1(\mathbb{Z}_p) = \mathbb{Z}_p^*.$$

For a closed point $x \in X/\mathbb{F}_p$, the image under $\rho_E$ of the geometric Frobenius conjugacy class $\mathrm{Frob}_x$ is given by

$$\rho_E(\mathrm{Frob}_x) = \alpha(x).$$

A theorem of Igusa says that the rank one representation $\rho_E$ is a surjective map onto $\mathbb{Z}_p^*$. As $x$ varies over $X$, the $p$-adic unit $\alpha(x)$ is given by the value of a rigid analytic function at the Teichmüller lifting of $x$. This rigid analytic function can be explicitly expressed in terms of the bounded solution of the Picard-Fuch differential equation of the family $f$. Such a bounded solution arising from a more general family of Calabi-Yau varieties is expected to contain important arithmetic information about the mirror map in mirror symmetry, see Lian-Yau [LY] for some positive results in this direction. We shall, however, not discuss this point of view here.

For an integer $k$, the above explicit description for $\rho_E$ shows that the L-function of the $k$-th tensor power of $\rho_E$ is given by

$$L(\rho_E^{\otimes k}, T) = \prod_{x \in X_0} \frac{1}{1 - \alpha(x)^k T^{\deg(x)}}, \tag{2.1}$$

where $X_0$ denotes the set of closed points of $X/\mathbb{F}_q$. Dwork's conjecture in this case is the following

THEOREM 2.1 ([D4]). *For each integer $k$, the L-function $L(\rho_E^{\otimes k}, T)$ is $p$-adic meromorphic.*

The key fact of the proof is that in this elliptic case, there is an excellent lifting as conjectured by Tate and proved by Deligne such that the rigid analytic function $\alpha(x)$ is overconvergent with respect to the excellent lifting, see [D3]. One then applies the Monsky trace formula to conclude the proof as noted by Dwork [D4].

Once we know that the L-function $L(\rho_E^{\otimes k}, T)$ is $p$-adic meromorphic. Many further arithmetic questions arise. For a connection between the special value $L(\rho_E^{\otimes k}, 1)$ and geometric Iwasawa theory, see Crew [Cr]. In the rest of this section, we shall discuss the $p$-adic absolute values of the zeroes of $L(\rho_E^{\otimes k}, T)$ and their arithmetic meaning in terms of modular forms.

For simplicity, our description will often be intuitive. An interested reader is encouraged to look up the references for precise definitions of some of the basic concepts. Let $k$ be an integer. Let $M_k$ be the space of overconvergent $p$-adic modular forms of level 2 and weight $k$. This is in general an infinite dimensional $p$-adic Banach space, see [K2] for a systematic treatment. It includes all classical modular forms of weight $k$ and level $2p^i$ for all integers $i > 0$. The Atkin $U_p$-operator acting on the Fourier expansion of a $p$-adic modular form is given by the map

$$U_p : \sum_{n \geq 0} a_n q^n \longrightarrow \sum_{n \geq 0} a_{np} q^n,$$

where $q$ (not a power of $p$) denotes the standard variable in the Fourier expansion of a modular form. The operator $U_p$ is a nuclear operator on $M_k$. Thus, the Fredholm determinant

$$D(k, T) = \det(I - U_p T | M_k)$$

is a $p$-adic entire function. Its coefficients are $p$-adic integers.

For arithmetic applications, it is important to understand the $p$-adic absolute values of the zeroes of the Fredholm determinant $D(k, T)$. One arithmetic reason comes about as follow. The Fredholm determinant $D(k, T)$ is the $p$-adic analogue of the Hecke polynomial $H_k(T)$ of the $p$-th Hecke operator associated to classical modular forms of weight $k$ and level 2. In fact, the Hecke polynomial $H_k(T)$ is essentially the part of $D(k, T)$ with small slopes, see Coleman [C1]. The Ramanujan-Peterson conjecture as proved by Deligne [De1] is to determine the complex absolute values of the zeroes of the Hecke polynomial $H_k(T)$. The $p$-adic version of the Ramanujan-Peterson conjecture is then to determine the $p$-adic absolute values of the zeroes of the Hecke polynomial $H_k(T)$. A good understanding of this $p$-adic question reveals subtle arithmetic information about the coefficients of modular forms. Unfortunately, the $p$-adic Ramanujan-Peterson conjecture seems to be much too hard in general. In fact, one does not even have a clean conjectural statement. The Gouvêa-Mazur conjecture discussed below can be viewed as a clean conjecture in this direction, not about a single Hecke polynomial $H_k(T)$ but about how the $p$-adic absolute values of the zeroes of $H_k(T)$ vary as the weight $k$ varies. One could work directly with $H_k(T)$ and use the $k$-th symmetric power of the first relative crystalline cohomology of the universal family of elliptic curves. We shall however focus our attention on $p$-adic modular forms and study the $p$-adic entire function $D(k, T)$ instead of the Hecke polynomial $H_k(T)$, because $D(k, T)$ contains more information than $H_k(T)$.

The relation between the L-function $L(\rho_E^{\otimes k}, T)$ and the Fredholm determinant $D(k, T)$ is given by the following theorem.

THEOREM 2.2. *For each integer $k$, we have the equality*

$$L(\rho_E^{\otimes k}, T) = \frac{D(k+2, T)}{D(k, pT)}. \tag{2.2}$$

*Equivalently,*

$$D(k, T) = \prod_{j \geq 0} L(\rho_E^{\otimes k-2-2j}, p^j T). \tag{2.3}$$

Equation (2.2) is really the Monsky trace formula applied to the rank one unit root F-crystal whose Frobenius matrix $\alpha(x)^k$ is overconvergent with respect to the canonical Deligne-Tate lifting, see [C2] and [K2] for a proof. The $U_p$-operator becomes the Dwork trace operator in this case. This gives a proof of Theorem 2.1. Equations (2.2)-(2.3) show that the L-functions $L(\rho_E^{\otimes k}, T)$ for all $k$ and the Fredholm determinants $D(k, T)$ for all $k$ determine each other. Thus, all results below can be formulated using either the L-functions $L(\rho_E^{\otimes k}, T)$ or the Fredholm determinants $D(k, T)$. We shall focus more on the Fredholm determinants $D(k, T)$ and indicate some of the translations to the L-functions $L(\rho_E^{\otimes k}, T)$.

One of Dwork's original motivations of proving Theorem 2.1 was to study $p$-adic properties of the Hecke polynomial $H_k(T)$, motivated in part by Ihara's work [Ih] relating Hecke polynomials to symmetric powers of elliptic curves. As explained above, this is related to arithmetic of modular forms such as some sort of $p$-adic Ramanujan-Peterson conjecture. Dwork, however, did not go any further in this direction. Motivated by a number of additional arithmetic applications such as congruences of modular forms, $p$-adic family of modular forms with a given slope and $p$-adic family of Galois representations of $\text{Gal}(\bar{\mathbb{Q}}/\mathbb{Q})$, Gouvêa-Mazur [GM1] proposed to understand how the $p$-adic absolute values of the zeroes of $D(k, T)$ vary as the integer weight $k$ varies $p$-adically.

The explicit definition of $L(\rho_E^{\otimes k}, T)$ in (2.1) and Fermat's little theorem show that $L(\rho_E^{\otimes k}, T)$ is $p$-adically continuous in $k$. More precisely, if $k_1$ and $k_2$ are two integers such that

$$k_1 \equiv k_2 \mod (p-1)p^m, \tag{2.4}$$

then

$$L(\rho_E^{\otimes k_1}, T) \equiv L(\rho_E^{\otimes k_2}, T) \mod p^{m+1}. \tag{2.5}$$

It follows from (2.3) that if $k_1$ and $k_2$ are two integers satisfying (2.4), then we also have

$$D(k_1, T) \equiv D(k_2, T) \mod p^{m+1}. \tag{2.6}$$

This gives a simple continuity result for $D(k, T)$ obtained independently by a number of authors, see [Ad], [Ko] and [GM2]. To get deeper information about the zeroes of $D(k, T)$, we would like to decompose $D(k, T)$ in terms of its slopes.

For a given integer $k$, the $p$-adic Weierstrass factorization theorem shows that the $p$-adic entire function $D(k, T)$ can be factored completely over $\Omega_p$:

$$D(k, T) = \prod_{i \geq 0} (1 - z_i(k)T), \tag{2.7}$$

where each $z_i(k)$ is a $p$-adic integer in $\Omega_p$. For a rational number $s \in \mathbb{Q}$, define the slope $s$ part of $D(k,T)$ to be

$$D_s(k,T) = \prod_{\mathrm{ord}_p z_i(k) = s} (1 - z_i(k)T). \tag{2.8}$$

This is actually a polynomial with coefficients in $\mathbb{Z}_p$. It is non-trivial only for $s \geq 0$. Let $d_s(k)$ denote the degree of the polynomial $D_s(k,T)$. This is a non-negative function in two variables. The function $d_s(k)$ is called the degree function of the entire function $D(k,T)$. The quantity $d_s(k)$ is the dimension of the space of overconvergent $p$-adic modular forms of weight $k$, level 2 and slope $s$. Thus, we may also call the degree function $d_s(k)$ as the dimension function of $p$-adic modular forms. A fundamental question is to understand the degree function $d_s(k)$, when it is positive, how it varies with $s$ and $k$, how large $d_s(k)$ can be.

Similarly, we can define another related degree function using the meromorphic L-function $L(\rho_E^{\otimes k}, T)$ instead of the entire $D(k,T)$. We now explain this point of view which will be used in our future work on a family of higher dimensional varieties where the automorphic form interpretation of Dwork's unit root L-function is not yet available. For a given integer $k$, the $p$-adic meromorphic function $L(\rho_E^{\otimes k}, T)$ can be factored completely over $\Omega_p$:

$$L(\rho_E^{\otimes k}, T) = \prod_{i \geq 0} (1 - u_i(k)T)^{\pm 1}, \tag{2.9}$$

where each $u_i(k)$ is a $p$-adic integer in $\Omega_p$. For a rational number $s \in \mathbb{Q}$, define the slope $s$ part of $L(\rho_E^{\otimes k}, T)$ to be

$$L_s(k,T) = \prod_{\mathrm{ord}_p u_i(k) = s} (1 - u_i(k)T)^{\pm 1}. \tag{2.10}$$

This is a rational function with coefficients in $\mathbb{Z}_p$. It is non-trivial only for $s \geq 0$. Let $d'_s(k)$ denote the degree of the rational function $L_s(k,T)$, which means the degree of its numerator minus the degree of its denominator. This function $d'_s(k)$ of two variables is called the degree function of the meromorphic L-function $L(\rho_E^{\otimes k}, T)$. Unlike $d_s(k)$, the new function $d'_s(k)$ can take negative values. We also want to understand $d'_s(k)$.

The relationship between the two degree functions is given by the following result, which is an immediate consequence of Theorem 2.2.

COROLLARY 2.3. *For each integer $k$ and each slope $s$, we have the equality*

$$d'_s(k) = d_s(k+2) - d_{s-1}(k). \tag{2.11}$$

*Equivalently,*

$$d_s(k) = \sum_{0 \leq j \leq s} d'_{s-j}(k - 2 - 2j), \tag{2.12}$$

*where $j$ runs over the rational numbers in the interval $[0, s]$.*

In order to understand how the two degree functions vary with $k$ for a fixed $s$, we need to introduce two more functions. For any rational number $s$, we define $m(E, s)$ to be the smallest non-negative integer $m \in \mathbb{Z}_{\geq 0} \cup \{+\infty\}$ such that whenever $k_1$ and $k_2$ are two integers satisfying (2.4), we have the equality

$$d_{s^*}(k_1) = d_{s^*}(k_2)$$

for all rational $s^*$ with $s^* \leq s$. Similarly, for any rational number $s$, we define $m'(E, s)$ to be the smallest non-negative integer $m \in \mathbb{Z}_{\geq 0} \cup \{+\infty\}$ such that whenever $k_1$ and $k_2$ are two integers satisfying (2.4), we have the equality

$$d'_{s^*}(k_1) = d'_{s^*}(k_2)$$

for all rational $s^*$ with $s^* \leq s$. By our definition, the following trivial inequalities

$$0 \leq m(E, s) \leq +\infty, \ 0 \leq m'(E, s) \leq +\infty$$

hold for each $s$. Furthermore, $m(E, s) = m'(E, s) = 0$ for $s < 0$. This is because $d_s(k) = d'_s(k) = 0$ for $s < 0$. The quantity $m(E, s)$ gives information about how often the degree function $d_s(k)$ changes as $k$ varies $p$-adically, where the slope $s$ is fixed. The quantity $m'(E, s)$ gives information about how often the degree function $d'_s(k)$ changes as $k$ varies $p$-adically, where the slope $s$ is fixed.

Using Corollary 2.3, one easily deduces the following comparison result.

COROLLARY 2.4. *For each integer $k$ and each slope $s$, we have*

$$m(E, s) = m'(E, s).$$

Thus, among the two quantities $m(E, s)$ and $m'(E, s)$, it suffices to understand one of them. We shall use $m(E, s)$ since the notation is simpler. Let $\lceil s \rceil$ denote the smallest integer which is at least as large as $s$. Motivated by Hida's work on a family of ordinary modular forms and based on some numerical computations, Gouvêa-Mazur [GM1] proposed the following rather precise conjecture about $m(E, s)$.

CONJECTURE 2.5 (GOUVÊA-MAZUR). *For all $s \in \mathbb{Q}_{\geq 0}$, we have*

$$m(E, s) \leq \lceil s \rceil.$$

This conjecture is true for $s = 0$ since $m(E, 0) = 0$ by the congruence in (2.6). A qualitative version of the conjecture was proved by Coleman [C2].

THEOREM 2.6 (COLEMAN). *For all $s \in \mathbb{Q}_{\geq 0}$, we have*

$$m(E, s) < +\infty.$$

This theorem shows that for each fixed slope $s$, the degree function $d_s(k)$ of $D(k, T)$ is a locally constant function of $k$ with respect to the $p$-adic topology of $\mathbb{Z}$ in each residue class modulo $p - 1$. More recently, it is shown in [W3] that the function $m(E, s)$ can be bounded by a quadratic polynomial in $s$. That is,

THEOREM 2.7. *There are two finite explicit constants $a_p$ and $b_p$ such that for all $s \in \mathbb{Q}_{\geq 0}$, we have the inequality*

$$m(E, s) \leq a_p s^2 + b_p s.$$

To prove this theorem, one first gives a uniform quadratic lower bound for the Newton polygon of $D(k, T)$. Then one transforms the uniform quadratic lower bound for the Newton polygon into a quadratic upper bound for $m(E, s)$ using a reciprocity lemma and the congruence (2.6).

I do not know if the quadratic bound in Theorem 2.7 could be improved. In fact, I do not even know if the dependence on $p$ of the constants $a_p$ and $b_p$ could be removed. Improvements about the quadratic bound would be possible if there is a sufficient amount of cancellation of zeroes on the right side of (2.2). If this is true, it would be more important to understand $L(\rho_E^{\otimes k}, T)$ than $D(k, T)$ since $L(\rho_E^{\otimes k}, T)$ would detect the cancellation of zeroes while $D(k, T)$ would not. From a heuristic cohomological point of view, $L(\rho_E^{\otimes k}, T)$ is an L-function which would have some sort of $p$-adic cohomological formula. In comparison, $D(k, T)$ is the characteristic series on the chain level and thus would contain redundant information about zeroes if there exists a cohomological formula. But this heuristic cohomological argument does not work since the involved horizontal connection of the unit root F-crystal is not overconvergent. This causes an essential difficulty in $p$-adic spectral theory. It explains why the chain level formula in (2.2) would not pass to any naively defined $p$-adic cohomology formula. Furthermore, the non-rationality of $L(\rho_E^{\otimes k}, T)$ would imply that there is a good portion of non-cancellation of zeroes on the right side of (2.2). Thus, the best one could hope for would be some sort of very partial $p$-adic cohomological formula which hopefully would explain some non-trivial cancellation of zeroes. We do not know how to proceed in this direction.

We shall instead propose a new conjecture concerning the slopes of the zeroes of $D(k, T)$, where the slope of a $p$-adic number $\alpha$ simply means $\mathrm{ord}_p \alpha$. This conjecture may be viewed as another step toward understanding how the degree function $d_s(k)$ varies as the slope $s$ and the weight $k$ vary. Let $S(k, p)$ be the set of slopes of the reciprocal zeroes of $D(k, T)$. This is in general an infinite set in $\mathbb{Q}_{\geq 0}$.

CONJECTURE 2.8 ([W3]). *For any given $p$, the set $S(k, p)$ has a uniformly bounded denominator for all $k$.*

This conjecture can be reformulated using the L-function $L(\rho_E^{\otimes k}, T)$. For this purpose, let $S'(k, p)$ be the set of slopes of the reciprocal zeroes and reciprocal poles of the L-function $L(\rho_E^{\otimes k}, T)$. This is in general an infinite set in $\mathbb{Q}_{\geq 0}$. Theorem 2.2 shows that we have the inclusion relations

$$S'(k, p) \subset S(k+2, p) \bigcup \{S(k, p) + 1\},$$

$$S(k, p) \subset \bigcup_{j \geq 0} \{S'(k - 2 - 2j, p) + j\}.$$

It follows that Conjecture 2.8 is equivalent to

CONJECTURE 2.9. *For any given $p$, the set $S'(k, p)$ has a uniformly bounded denominator for all $k$.*

Conjecture 2.8 easily translates into a clean non-trivial general statement about the $p$-adic absolute values of the reciprocal zeroes of all the Hecke polynomials $H_k(T)$. Thus, it can be viewed as a $p$-adic Ramanujan-Peterson conjecture in a suitable sense. Except for the trivial case when $D(k, T)$ is already a polynomial, I do not know a single non-trivial example for which Conjecture 2.8 is true, even for a fixed $k$. On the other hand, the classical Riemann hypothesis says that the real part of the zeroes of the Riemann zeta function is a rational number whose denominator is bounded by 2. The real part of a zero for the Riemann zeta function corresponds exactly to the slope of a reciprocal zero in our $p$-adic situation. Thus, Conjecture 2.8 can also be viewed as a $p$-adic Riemann hypothesis for the elliptic family. We do not conjecture an explicit bound for the denominator of $S(k,p)$. Any proof or a sufficient amount of numerical computations of Conjecture 2.8 would likely produce such an explicit bound.

Finally, we turn to the size problem of the degree function $d_s(k)$. We want to discuss possible uniform finiteness for the values of $d_s(k)$ and how this finiteness relates to Conjecture 2.8. As a preliminary evidence, one easily deduces from Theorem 2.6 that for bounded $s$, the function $d_s(k)$ is a bounded function of $k$. We raise the following much stronger

QUESTION 2.10. *Is the degree function $d_s(k)$ uniformly bounded for all $s$ and all $k$?*

For positive $d_s(k)$, the integer $d_s(k)$ is a denominator (not necessarily the smallest one) for the slope $s$ of the reciprocal zeroes of $D_s(k,T)$. Thus, a positive answer of Question 2.10 implies that the set $S(k,p)$ has a uniformly bounded denominator for all $k$, which is precisely what Conjecture 2.8 says. We are inclined to believe that Question 2.10 has a positive answer in the current elliptic family case, although we do not believe it in general higher dimensional case. To give some evidence why Question 2.10 might have a positive answer, we include the following simple result which shows that on the average, Question 2.10 already has a positive answer.

PROPOSITION 2.11. *There is an explicit positive constant $c_p$ such that for all real numbers $A \geq 1$, the inequality*

$$\frac{1}{A} \sum_{0 \leq s \leq A} d_s(k) \leq c_p$$

*holds uniformly for all $k$, where $s$ runs over $[0,A] \cap \mathbb{Q}_{\geq 0}$.*

**Proof**. By [W3], there is an explicit uniform quadratic lower bound for the Newton polygon of $D(k,T)$. This implies that there are two explicit positive constants $c_p > 0$ and $e_p$ such that

$$\sum_{0 \leq s \leq A} s d_s(k) \geq \frac{1}{c_p} (\sum_{0 \leq s \leq A} d_s(k))^2 - e_p (\sum_{0 \leq s \leq A} d_s(k)). \qquad (2.13)$$

If

$$\sum_{0 \leq s \leq A} d_s(k) = 0,$$

Proposition 2.11 is trivially true. If

$$\sum_{0 \leq s \leq A} d_s(k) \neq 0,$$

we can cancel it from (2.13). Using the condition $0 \leq s \leq A$, we deduce

$$A \geq \frac{1}{c_p} (\sum_{0 \leq s \leq A} d_s(k)) - e_p. \qquad (2.14)$$

This implies that

$$\frac{1}{A} \sum_{0 \leq s \leq A} d_s(k) \leq c_p + \frac{c_p e_p}{A}.$$

Increasing the size of $c_p$ if necessary, we conclude that Proposition 2.11 is true.

The above result shows that on the average, there are at most $c_p$ reciprocal zeroes of $D(k,T)$ whose slopes are in a given interval $I \subset \mathbb{R}_{\geq 0}$ of length 1. It also shows that for any given $k$, most of the non-negative integers $d_s(k)$ for $s \in [0, A]$ are bounded by $c_p$. Taking $A = s$ in Proposition 2.11, we deduce the following weak but simple

COROLLARY 2.12. *There is an explicit positive constant $c_p$ such that for all $s \geq 0$ and all integers $k$, we have the uniform linear bound*

$$d_s(k) \leq c_p s.$$

## 3. Some interactions with B. Dwork

This section describes some of my personal interactions with B. Dwork and how some of his mathematical work had influenced my study of mathematics.

My initial interest in Dwork's $p$-adic theory grew out of my attempt to understand diophantine equations over a finite field. When I was still a graduate student in mid-eighties working under the direction of Neal Koblitz at the University of Washington, I became fascinated with the simple but beautiful theorem of Chevalley-Warning, which counts the number of rational points modulo the characteristic $p$. In order to get the full information about the solution number and to see how it varies in various extension fields as predicted by the Weil conjectures, it was very natural to try to lift the argument in a systematic way to characteristic zero. I spent several months trying to lift the Chevalley-Warning argument but could not control it in a systematic way when the base field varies. Then I realized that the difficult I had was already succeesfully overcome by Dwork [D1] before I was born. The systematic lifting led Dwork to his fundamental trace formula, which is the key toward his rationality proof of the zeta function of an algebraic variety over a finite field. Years later, I mentioned this in a conversation with Dwork. He told me that he was indeed mostly influenced by Warning in his rationality proof.

The first time I met Dwork was at the AMS 1989 meeting in Muncie, Indiana. By that time, I already knew his trace formula in the classical overconvergent setting and was able to use it to obtain some preliminary information about the Adolphson-Sperber conjecture [AS] on the generic Newton polygon for the L-function of exponential sums. This led to my later result [W1] which proves a modified form of the Adolphson-Sperber conjecture and gives a systematic way to determine when the generic Newton polygon coincides with its lower bound (the Hodge polygon). In particular, it provided a direct $p$-adic proof of Mazur's conjecture [Ma] which says that the Newton polygon coincides with the Hodge polygon for a generic hypersurface.

Before the Muncie meeting, I had also made some natural experimental study myself about meromorphic continuation of what I called formal L-functions. My intention of such a formal study was to test and to see how far one can go with Dwork's trace formula. I showed an optimal $c\log$-convergent result for the formal L-function I defined. In particular, a counter-example was found which shows the formal L-function is not always a meromorphic function. I explained my results to Dwork at the meeting. He said that he would believe them. It turned out that he and Sperber [DS] had also proved a similar $c\log$-convergent result for formal L-functions, although they could not prove that their result is optimal in general, perhaps partly due to Katz's more general meromorphic conjecture about the L-function of an F-crystal.

I felt that my formal counter-example could already be a counter-example for Katz's meromorphic conjecture. But at that time, I did not understand the rather fancy definition of F-crystals and thus I could not check that the formal L-function I studied essentially agrees with the L-function of an F-crystal. Several years later, in 1993 at the Igusa retirement conference held at the Johns Hopkins University, I met Katz and asked him to explain the concrete meaning of an F-crystal. According to what Katz explained, I was more convinced that my counter-example should be a counte-example for his meromorphy conjecture about the L-function of an F-crystal. I mentioned this to Dwork again at the Igusa retirement conference (Katz already left the conference and returned to Princeton). This time, he encouraged me to write it up and send it to Katz. The interest of this work [W2] to me was to see the limit of the Dwork trace formula and to know where to stop.

At the Muncie meeting, I also learned Dwork's meromorphic conjecture for the L-function of his unit root F-crystal (a continuous $p$-adic representation) arising from algebraic geometry. This conjecture was proposed in early seventies in his attempt to understand how the roots of the zeta function of a family of varieties vary when the parameter varies. The ordinary Hodge-Newton decomposition as he proved in [D5-6] connects the slopes of the Frobenius eigenvalues with $p$-adic reducibility of the associated Picard-Fuch differential equation. This line of investigation was continued and essentially completed by Katz [K3]. In some cases, Dwork could even give an explicit analytic formula for the unit root part of the zeta function purely in terms of the $p$-adically bounded solutions of the Picard-Fuch equation. He attached an L-function to such a unit root formula (which is a continuous $p$-adic representation of some $\pi_1^{\mathrm{arith}}(X)$) and conjectured its $p$-adic meromorphic continuation. As discussed in section 1, one version of the conjecture simply says that the L-function of the relative $p$-adic ètale cohomology of a family of varieties over a finite field of characteristic $p$ is a $p$-adic meromorphic function. From family point of view, Dwork's conjecture can be viewed as the starting point of a truly $p$-adic extension of the Weil conjectures from a family of zero dimensional varieties to a family of positive dimensional varieties.

Such a geometric unit root L-function seems quite mysterious since the Frobenius map of the unit root F-crystal is no longer overconvergent. From formal L-function point of view, the geometric unit root L-function looks as bad as a general formal L-function. Thus, my counter-example for the Katz conjecture shows that there is nothing more one can say along the direction of formal L-functions. In several interesting geometric cases treated by Dwork, such as the family of ordinary elliptic curves as described in section 2, the universal family of ordinary genus 3 plane curves and a certain family of ordinary K-3 surfaces [D6], he was able to get around the difficulty by showing that the involved unit root F-crystal (its Frobenius map) is in fact overconvergent with respect to another lifting called excellent lifting. Thus, in such a case, the situation is reduced to the classical overconvergent case. But, unfortunately, excellent lifting rarely exists and thus the conjecture cannot be reduced to the "trivial" overconvergent situation in general. Even in the few exceptional cases where the excellent lifting does exist, the situation is quite subtle as one might guess from the simplest example discussed in section 2.

Right after the Muncie meeting, I started to have some preliminary feeling about Dwork's conjecture. I had a one-line argument which already proves something slightly stronger than the result in [DS]. However, for a long time I could not see why Dwork's conjecture should be true in the general case. Both the formal L-function approach and Dwork's excellent lifting approach looked hopeless to me. In the latter case I must say that I did not and still do not understand the concept (excellent lifting) very well. This could be one of the reasons that I decided to look for other approaches. During 93-94 academic year, I was a member at the Institute For Advanced Study at Princeton and had the chance to discuss some of my ideas with Deligne and Katz. The discussions were very helpful. A few months later, I felt that I had enough ideas to prove some essential rank one case of Dwork's conjecture such as the higher dimensional Kloosterman sum family. However, I still did not have a good feeling about the higher rank case of Dwork's conjecture and I was not sure if there would be a sufficient amount of interest in the partial rank one case that I could prove. For these and other reasons, I did not even try to write down my rank one ideas. In fact, for a while, I shifted my interests to other problems since I had no ideas how long it would take me to get the higher rank case. However, the problem was on my mind and I knew that I would return to it earlier or later.

During 93-94, I also had a chance to talk to Dwork. He had an appointment in Italy but returned to Princeton for a short period. I explained to him some of the possibilities of proving his conjecture. He did not have a good feeling about my possible approach and he did not ask me to explain in further detail. Instead, he suggested to me to look at the family of ordinary abelian varieties, generalizing the elliptic family. He had strong intuition to feel that the excellent lifting should exist in the abelian variety case but he did not know enough about abelian varieties. He said that the proof in this case would be a respectable work although he would not consider it to be great. I was not familiar with abelian varieties either. Even worse, I was not even familiar with excellent lifting. I was pretty sure that I would not be able to do anything with it. Thus, I did not take his suggestion seriously. In the end, he gave me a pack of his old reprints. I said "that would keep me busy for the rest of my life". He said "just for a few days".

By 1996, the Gouvêa-Mazur conjecture and Coleman's work about modular forms caught my attention because they are closely related to the unit root L-function studied by Dwork. Although Dwork's conjecture was already known in the elliptic family case, the connection with the Gouvêa-Mazur conjecture suggested another hint about the potential significance of Dwork's conjecture in general. This gave me additional motivations to return to the investigation of Dwork's conjecture. During the summer of 97 when I was visiting Sichuan University in Chengdu and the Mathematics Institute in Beijing, I had a chance to return to the problem and found additional ideas which allowed me to handle the whole rank one case of Dwork's conjecture. But I still did not have a good feeling about the higher rank case.

At the 97 Seattle summer research conference on finite fields, I spoke about Dwork's conjecture as described in the previous two sections. I did not announce that I would be able to prove the rank one case of Dwork's conjecture, partly because I ran out of lecture time and partly because I did not write down my ideas yet. After the conference, I started to work out the details of my rank one ideas and write up the rank one proof to make sure it is correct. After finished the first draft of the rank one case in a simpler setup, the global picture was coming together. After a few more weeks of intensive thinking, I was starting to get some feeling on higher rank case as well. By the end of 97, I was confident that I had found all the main ideas that are needed for the proof of Dwork's conjecture in the general case. Although it would take a couple of more years to fill in all the details of the general proof, Dwork was quite excited to hear about it and wanted to understand the proof himself, partly due to his curiosity about how I could avoid his excellent lifting completely.

I was particularly moved that despite his serious illness, Dwork still came to visit Irvine in late January 1998, gave two excellent talks and had a lot of interesting discussions with me about his conjecture. He felt that the proof is correct but indicated that the detail has to be checked. For him, he mostly wanted to understand the rank one proof since the rank one case was already quite essential to him. He did not bother to try to understand the reduction from higher rank case to rank one case. He was aware that his subject area is not a popular one. Thus, he asked me why I was interested in proving his conjecture even though he did not ask me to do so. I told him that I like the conjecture and I think it is a good problem.

The last time I saw Dwork was in late March 1998 in his temporary Menlo Park home in San Francisco area. I stayed in their house for two days and enjoyed more discussions with him during my stay. His illness was more serious. He would more easily fall to sleep for a few minutes from times to times. When he woke up, he would continue his conversation from where he stopped. My impression was that he still had a very strong intuition and a young curious mind. It was still intellectually interesting and rewarding to talk to him. If I mention a problem that is interesting to him, he would try to think about it and suggest what his feeling would be. There was one question which was bothering him for some time and which he thought to be a gap in his 1966 paper [D2] on the finite dimensionality of his cohomology space for a singular hypersurface. He explained the problem to me. Then, I explained to him why it did not seem to be a gap to me. I mentioned that very likely he used the same argument more than thirty years ago and it was too easy for him to write it down. He felt a little better but was not fully convinced. A week later, he sent me an email happily saying that he understood the problem and had fixed it. I was very pleased in knowing that the problem was no-longer bothering him.

In early April 98, Dwork and his wife Shirely moved back to Princeton. I was hoping to visit him again in September 98 if he would be well enough. For a while, I did not hear from him. On May 8, 1998, I received an email message written by his son saying that Dwork would be very happy to see me in September if he would be well by that time. I was very pleased. Two days later, I was shocked to hear the news that Dwork died on May 9, 1998.

Looking back, I felt very fortunate to be able to have many discussions with Dwork during the last few months of his life. I was hoping that he would be able to recover from his illness so that some day he could explain to me some of his other results and conjectures which are probably not well known and could be quite difficult to read on one's own. That day would now never come. On the other hand, just by looking at the small fraction which I do know, it seems certain to me that Dwork's work would be inspiring for many years to come. It may take much longer to fully appreciate some of his insights.

Department of Mathematics, University of California, Irvine, CA 92697-3875
*E-mail address*: dwan@math.uci.edu